\documentclass[12pt,leqno]{amsart}

\def\opn#1#2{\def#1{\mathop{\kern0pt\fam0#2}\nolimits}}

\opn\arg{arg} \opn\DD{\mathcal D} \opn\GCD{GCD}
\opn\Id{Id} \opn\im{Im} \opn\rad{rad}
\opn\sgdeg{sgdeg}  \opn\Frac{Frac} \opn\Spec{Spec}

\def\:{:\text{ }}

\begin{document}

\begin{center}{\large \bf Adjamagbo Determinant and Serre Conjecture
\\[.2cm]
for linear groups over Weyl algebras}\\[0.2cm]
{Kossivi Adjamagbo \\ Universit\'e Paris 6 - Case 172 - Institut de Mathématiques de Jussieu\\ 4, place Jussieu, 75252 PARIS CEDEX 05\\ adja@math.jussieu.fr}
\end{center}

\vspace*{.2cm}

{\tiny {\bf Abstract:} Thanks to the theory of determinants over an Ore domain, also called Adjamagbo determinant by the Russian school of non commutative algebra, we extend to any Weyl algebra over a field of characteristic zero Suslin theorem solving what Suslin himself called the $K_1$-analogue of the well-known Serre Conjecture and asserting that for any integer $n$ greater than 2, any $n$ by $n$ matrix with coefficients in any algebra of polynomials over a field and with determinant one is the product of elementary matrices with coefficients in this algebra .}

\vspace*{.5cm}

\noindent
{\bf Introduction}

Let $A$ be a ring, $n$ a positive integer, $M_n(A)$ the ring of $n$ by $n$
matrices with
coeficients in $A$, $GL_n(A)$ the group of invertible elements of $M_n(A)$
and $E_n(A)$
its subgroup generated by its elementary matrices.

Let us consider the following natural question : Does there exist a
function of the
coefficients of the elements of $M_n(A)$ with values somewhere which allows the
characterization of elements of $E_n(A)$ between those of $M_n(A)$ and not
only between those
of $GL_n(A)$ and such that this function is ``effectively computable''
whenever the ring $A$
is ``effective''  in the sense of \cite{adja3}?

Let us first assume $A$ commutative. In this case, we think naturally to
the determinant
function of matrices in $M_n(A)$. Indeed, for any element $a$ of $M_n(A)$,
if $a \in E_n(A)$,
then it is obvious that $\det a = 1$. But the converse is false in general,
as shown by the
famous counter-example of P.M. Cohn \cite{cohn}, prop. 7.3 :
$$ a = \begin{pmatrix} 1 - XY & - Y^2 \\ X^2 & 1 + XY \end{pmatrix}\in
M_2({\Bbb Q}[ X,Y ])$$
Nevertheless, it is well known that this converse is true if $A$ is a
field, or an euclidian
ring or a semi-local one, or if $A$ is a noetherian ring whose Krull
dimension is lower than
$n - 1$, according to Bass-Milnor-Serre theorem on products of elementary
matrices published in 1967 \cite{bms}. 

On the other hand, Suslin proved ten years later in \cite{suslin} that
this last thoerem can
be improved by taking $n$ only greater than 2 (instead of 1 plus the Krull
dimension of $A$)
and by choosing $A$ as any algebra of polynomials over a field. In the introduction of this paper, Suslin himself presented this result as the solution to ``the $K_1$-analogue of the well-known serre Problem recently solved completely by the author, and independently by Daniel Quillen''. Indeed, according to the triviality of the special Withehead group $SK_1 A$ of such algebra, for a square matrix $a$ with entries in $A$, $\det a = 1$ means that this matrix is ``stably'' (i.e. after augmentation of the matrix by adding some 1's on the diagonal and 0 outside) a product of elementary square matrices with entries in $A$, and the problem is to know if it is ``actually'' a product of elementary square matrices with entries in $A$. This justifies Suslin's analogy with Serre problem which ask if any ``stably'' (i.e. after addition by a finite free $A$-module) free A-module of finite type is ``actually'' free.

Thanks to the theory of determinants over an Ore domain developped in \cite{adja4}, summed up in
\cite{adja5}, more brievely in \cite{bjork}, A,III,  and already called ``Adjamagbo determinant'' by Russian school of non commutative algebra following A. Mikhalev and A. Guterman (see for instance \cite{guntmik1}, \cite{guntmik2}, \cite{gunter1}, \cite{gunter2}, \cite{gunter3}), we extended in \cite{adja5},  this theorem of Bass-Milnor-Serre to
the case where $A$ is
a non-commutative ``classical filtered ring'', i.e. a ring endowed with an
increasing ${\Bbb N}$
- filtration $F$ whose associated graded ring is a commutative regular
domain flat over its
subring F(0) which has a trivial special Whitehead group $SK_1 F(0)$. It is
in particular the
case of the envolopping algebra of a Lie algebra of finite dimension over a
field. It is also
the case of a classical or a formal (resp. analytic) Weyl algebra over a
field (resp. the field
of real or complex numbers), see for instance \cite{adja5}, p. 404.

Exactly as Suslin did for Bass-Milnor-Serre theorem, the aim of this note
is to improve this
generalization of Bass-Milnor-Serre theorem, by taking $n$ only greater
than 2 (instead of 1
plus the Krull dimension of $A$) and by choosing $A$ as any Weyl algebra
over a field of
characteristic zero, thanks to Stafford Main Theorem on the `` module
structure of
Weyl algebras'' \cite{stafford} and to Varerstein $K_1 $-stability theorem \cite{vaserstein}, th. 3.2.

Finally, we formulate a natural conjecture about filtered rings which would
prove that this Suslin
Theorem follows from the forcomming theorem. We end with other open
problems related to this theorem.
\\ \\
{\bf Recall} (on the canonical determinant over filtered Weyl algebras, see for instance \cite{adja5})

1) $K$ being a commutative field and $m$ a positive integer, let us denote by
 $(X_i, Y_i)_{1\leq i \leq m}$ a system of indeterminates over $K$,
 $K[X_1, \ldots, X_m]$ the $K$-algebra of polynomials in indeterminates
$X_1, \ldots, X_m$,
 $A_m (K)$ the Weyl algebra over $K$ of index $m$,
 i.e. the $K$-algebra $K[X_1, \ldots, X_m][\partial / \partial X_1, \ldots,
\partial / \partial
X_m]$ of $K$-linear differential operators over $K[X_1, \ldots, X_m]$, $F$ the
differential filtration on $A_m (K)$, i.e. $F(0) = K[X_1, \ldots, X_m]$ and
$F(j + 1) = F(j) +
\sum_{1 \leq i \leq m} F(j) \partial / \partial X_i$ for each natural
integer $j$, $gr_F A_m K$
the associated graded ring, $gr_F$ the canonical map from $A_m (K)$ to
$gr_F A_m K$ and $gr_F
(A_m (K))$ the image of $gr_F$, i.e. the set of homogeneous elements of the
graded ring $gr_F
A_m K$, or in terms of partial differential equations, the set of principal
symbols of the
differential operators belonging to $A_m (K)$.

2) The ring $gr_F A_m (K)$ being isomorphic to the regular ring of
polynomials over $K$ in
indeterminates $X_1, \ldots, X_m, Y_1, \ldots, Y_m$ and the group $SK_1
F(0) = SK_1 K[X_1,
\ldots, X_m]$, isomorphic to $SK_1 K$ according to Quillen Theorem \cite{quillen}, §
6, th. 7, being
trivial, $F$ is therefore a ``classic regular filtration'' on
$A_m (K)$ making the later be a ``classic regular ring'' and in particular an
Ore domain.

3) According to the theory of determinants over Ore domains, there exists
an unique map from
the set of square matrices of elements of $A_m (K)$ and with values in
$gr_F (A_m (K))$, denoted
by ${\det}_F$ and called ``the canonical determinant over the filtered Weyl
algebra over $K$
of index $m$'' or ``the principal determinant over the Weyl algebra over
$K$ of index $m$'',
such that, for any square matrices $a$, $b$ of the same size with
coefficients in
$A_m (K)$ and any diagonal matrix $diag( x, 1,\ldots,1)$ with diagonal
coefficients
$x, 1,\ldots, 1$ in $A_m (K)$, we have:
$$ {\det}_F (a b) = {\det}_F (a) {\det}_F (b) \; \; (the \; homomorphism \;
axiom)$$
$$ {\det}_F ( diag( x, 1,\ldots, 1)) = gr_F (x) \; \;  (the \; prolongation
\; axiom)$$

4) The first fondamental property following from this homomorphism axiom is
that for
any elementary $e$ with coefficients in $A_m (K)$, we have:
$${\det}_F (e) = 1 \; \; (the \; elementary \; property)$$

5) The second fondamental property following easily from these two axioms
and this elementary
property is that for any triangular matrix $t$ with coefficient in $A_m
(K)$, ${\det}_F (t)$ is
the product of the principal symbols of its diagonal coefficients (the
triangular property)

6) The third fondamantal property following from these two axioms and this
triangular property
is that for any $n$ by $n$ matrix $a$ with coefficients in $A_m (K)$,
${\det}_F (a)$ can be
computed in a pratical way by ``Gauss method''. Indeed, thanks to the left
commun multiple
property of the Ore domain $A_m (K)$ and to suitable combinations on the
lines of square
matrices with coefficients in $A_m (K)$, it is possible to find $n$ by $n$
elementary
or cancellable diagonal matrices $p_1, \ldots, p_r$ with coefficients in
$A_m (K)$ such that
$p_1 \ldots p_r a$ is an upper triangular matrix $t$. Then, thanks to the
homomorphism axiom
and the triangular property, we obtain the following explicit expression of
${\det}_F (a)$ as
quotient of two principal symbols of elements of $A_m (K)$ which can
finally be reduced to the
principal symbol of a element of $A_m (K)$ thanks to the factoriality of
the ring $gr_F A_m(K)$ and to the ``regularity theorem'' :
\[{\det}_F (a) = {\prod}_{1 \leq i \leq n} gr_F (t(i,i)) / {\prod}_{1 \leq
k \leq r}
{\prod}_{1 \leq i \leq n} gr_F (p_k (i,i))\]

7) If the field $K$ is ``effective'' in the sense of \cite{adja2}, then it follows
from this last
formula and from the ``effectivity'' of the Ore property of the
``effective ring'' $A_m (K)$
that the restriction of ${\det}_F$ to each ``effective ring'' of square
matrices with
coefficients in $A_m (K)$ is an ``effectively computable'' function, as
proved in \cite{adja1}.

8) It also follows from this formula that if $a$ and $b$ are two square
matrices with
coefficients in $A_m (K)$ and if
$$ a \oplus b = \begin{pmatrix} a & 0 \\ 0 & b \end{pmatrix}$$
then $$ {\det}_F (a \oplus b) = {\det}_F (a) {\det}_F (b)$$

9) The link between the classical determinant over a commutative ring and
${\det}_F$ is that if
$B$ is any commutative sub-ring of $A_m (K)$ and $\det$ the classical
determinant over $B$,
then :
$${\det}_F {\mid}_B = gr_F \circ \det$$

10) Finally, one of the most remarkable analogies already discovered
between the classical
determinant over a commutative ring and $\det_F$ is that an element $a$ of
a ring of square
matrices over $A_m (K)$ is invertible if and only if ${\det}_F (a)$ is
invertible in $gr_F A_m(K)$, i. e. a non zero element of $K$.
\\ \\
{\bf Theorem}

For any integer $n$ greater than 2, a $n$ by $n$ matrix $a$ with
coefficients in any Weyl
algebra over any field of characteristic zero is a product of elementary matrices with coefficient in this algebra
if and only its
canonical determinant over this algebra is 1 .\\ \\
{\bf Proof}

1) Let $A_m (K)$ be such a Weyl algebra and $F$ its differential filtration.

2) If $a$ is a product of elementary matrices with coefficients in $A_m
(K)$, then according to
the homomorphism axiom and the elementary property of ${\det}_F$, it is
obvious that ${\det}_F
(a) = 1$.

3) So let us assume conversly that ${\det}_F (a) = 1$. According the above
characterization of
elements of $GL_n (A_m (K))$ between those of $M_n (A_m (K))$, $a \in GL_n
(A_m (K))$.

4) $F$ being a classical regular filtration as we remarked above, it
follows from the cited
Quillen theorem that the canonical map from $K_1 F(0)$ in $K_1 A_m (K)$ is
an isomorphism. So
there exists a positive integer $p$, a $p$ by $p$ matrix $b$ with
coefficients in $F(0)$
and integers $r$ and $s$ such that $n + r = p + s$ and :
$$ (*) \; \; (a \oplus i_r)^{-1} (b \oplus i_s) \in E_{n + r} (A_m (K))$$
where $i_j$ denotes the unit of the group $GL_j (A_m (K))$ .

5) According to (*) and to the points 3), 4), 8) and 9) of the recall, we
have :
$${\det}_F (b \oplus i_s) = \det (b \oplus i_s) = 1$$

6) Since $F$ is a classical regular filtration, in particular since $SK_1
F(0)$ is trivial, it
follows from this last equality that there exists a positif interger $t$
such that :
$$ b \oplus i_{s + t} \in E_{p + s + t} (F(0))$$

7) So according to (*), we have :
$$ a \oplus i_{r + t} \in E_{n + r + t} (A_m (K))$$

8) On the other hand, according to Stafford cited theorem, the stable range
of $A_m (K)$ is 2.
Thanks to Vaserstein cited theorem, it follows from the last relation as
desired that :
$$a \in E_n (A_m (K))$$
Q.E.D.
\\ \\
\noindent
{\bf Remark 1}

1) The previous theorem may by interpreted in terms of systems of partial differential equations with solutions in any $K[X_1, \ldots, X_m]$-algebra $B$ which is a $A_m (K)$-left-module, more precisely in terms of such a system which is ``elementary resoluble by derivation'', following and refining \cite{adja1}.

2) The statement of the previous theorem is clearly the faithfull
generalization to Weyl algebras over
a field of characteristic zero of Suslin theorem for algebras of
polynomials over such a field.

3) Furthermore, it solves the $K_1$-analogue of Serre Conjecture over Weyl algebras, since according to points (2) and (7) of the previous proof, for any  square matrix $a$ with entries in a Weyl algebra over any field of characteristic zero $A$, ${\det}_F (a) = 1$ means that this matrix is ``stably'' a product of elementary square matrices with entries in $A$.   

4) On the other hand, Suslin theorem would follow from the previous theorem, with a non
commutative algebraic
proof completely different from Suslin original commutative algebraic
proof, if the following
``natural'' conjecture could be confirmed :
\\ \\
{\bf Conjecture}

If $A$ is a ring endowed with an increasing ${\Bbb N}$
- filtration $F$ such that the associated graded ring is a domain, then for
any integer $n \ge 2$, the
sub-ring $F(0)$ of $A$ verifies : $$ GL_n(F(0)) \cap E_n(A) = E_n(F(0))$$
\\ \\
\noindent
{\bf Remark 2}

1) The assumption on the associated graded ring means that the associated
degree function $deg_F$ is
"additive", i.e. is an homomorphism from $(A-\{ 0 \}, \times)$ to $({\Bbb
N}, +)$. Morever,
this degree function is such that $F(0) - \{ 0 \}  = {deg_F}^{-1}(0)$.

2) Using this degree fonction, it seems that the proof of the conjecture
could be purely formal, as
this statement could be easily checked for $n=2$ in the case where the
considered element of
$GL_n(F(0)) \cap E_n(A)$ is a product of at most five elements of $A$.

3) The main interest of this conjecture is that thanks to it, a non trivial
property of a
commutative ring (Suslin theorem) could be deduced from the similar
property of a ``simple'' non
commutative extension of this ring (previous theorem), in an analogous way
as some deep properties of
the field of real numbers could be deduced from the similar ones of an
``algebraically closed''
extension of this field like the field of complex numbers.

4) This kind of contribution of non commutative algebra to commutative
algebra seems non common in
mathematics. A confirmed example of such a contribution is the fact that
the famous Jacobian
Conjecture, which claim that any endomorphism of an algebra of polynomials
over a field of
characteristic zero with a non zero jacobian in this field is an
automorphism, could be deduced from
Dixmier Conjecture, which claims that any endomorphism of a Weyl algebra
over a field of
characteristic zero is an automorphism (see for instance \cite{bcw}, p. 297).

5) Another interest of the proposed conjecture is that it would prove that
Cohn counter-example
cited in the introduction is even better than one thinks now, in the sense
that it is not in
$E_2(A_2({\Bbb Q}))$, showing in this way that the lower bound 3 for $n$ in
the previous theorem is
the finest as in Suslin theorem.

6) Conformly to the ``effectiveness'' problem evocated in the introduction
and the recall, a natural
question risen from the previous theorem is the following :

For an integer $n$ greater than 2, how to split ``effectively'' a $n$ by
$n$ matrix  with
coefficients  in a Weyl algebra over an ``effective'' field of
characteristic zero and with principal
determinant 1 as a product of elementary ones ?

7) According to the prominent part that the cited Stafford stable rank
theorem plays in the proof of
the previous theorem, it is easy to conjecture that this last
``effectiveness'' problem should need the
resolution of the following one :

Given three elements $a, b, c$ of a Weyl algebra over an ``effective''
field, how to
find ``effectively'' two elements $d$ and $e$ of this algebra such that $a
+ d c$ and $b + e c$
generate that same left ideal of this algebra as $a$, $b$ and $c$ ?

8) Let us now consider a last question risen from the previous theorem.

Indeed, for the problems
met in the algebraic theory of partial differential equations, the working
noetherian domains of
differential operator are not Weyl algebras, but what could be called
``formal (resp. convergent)''
Weyl algebras, deduced from ``classical'' Weyl algebras by replacing
polynomials by formal (resp.
convergent) power series (see for instance \cite{adja5}, p. 404). Since it is
natural to hope that the
previous theorem works also for these ``working'' Weyl algebras, the proof
of this theorem lead us to
the following question :

Is the stable rank of a ``formal'' (resp. ``convergent'') Weyl algebra over
a field of
characteristic zero (resp. a sub-field of the field of complex numbers)
also 2 ?
\\ \\

\end{document}